\documentstyle[11pt,amsmath,amsfonts,here ]{article}

\newtheorem{th}{Theorem}

\def\R{{\mathbb R}}

\def\Z{{\mathbb Z}}

\def\a{\alpha}
\def\b{\beta}

\def\D{\Delta}
\def\rw{\rightarrow }

\def\d{\delta}

\def\s{\sigma}
\def\n{\newline}

\def\M{{\cal M}}

\def\ov{\overline}
\def\wh{\widehat}
\def\wt{\widetilde}
\def\wf{\widehat{f}}
\def\fr{ \ fr}
\def\bef{$\emptyset$-definable basic $f$-set}
\def\ef{$\emptyset$-definable $f$-set}

\begin{document}
\title{\large \bf ON O-MINIMALITY OF EXTENSIONS OF 
$\R$ BY RESTRICTED GENERIC SMOOTH FUNCTIONS.}
  \author{\large Alexei Grigoriev}
  \date{ 26.4.2005}
  \maketitle
  \begin{abstract}
It is shown that the extension of $\R$ by a generic smooth function restricted to the unit cube is o-minimal.
The generalization to countably many generic smooth functions is indicated. Possible applications are sketched.   
\end{abstract}
  
\bigskip
\bigskip
\bigskip
{\bf 0. Introduction.}\n

A key result in the theory of subanalytic sets is Gabrielov's
theorem ([Gab]), which states that the complement of a subanalytic set is again subanalytic.
One form to state this fact in the framework of mathematical logic, is as follows. 
Denote by ${\cal L}$ the language $\{0,1,+,-,\cdot,<\}$. For each function $f$ which is real analytic
on a neighbourhood of the closed unit cube of $\R^n$ for some $n$, add to ${\cal L}$ the symbol $\hat{f}$.
Denote the language obtained by $L_{an}$, and denote by $\R_{an}$ the set of the real numbers seen as a
$L_{an}$ structure, with each of the symbols $\hat{f}$ being interpreted as a function which is
equal to $f$ on the unit cube and to $0$ outside. Denote by $T_{an}$ the theory of $\R_{an}$.
Gabrielov theorem implies \n

{\bf Theorem [vdD1]. }{\it 
 $T_{an}$ is model complete and o-minimal.\n
}

Our purpose here is to prove an analogue of this theorem for extensions of $\R$  
by (countably many) generic {\it smooth} functions restricted to the unit cube.
To avoid overly complex notations, the argument is written down in detail for the case of extension 
by a single generic smooth function. We then comment on the generalization to countably many 
generic smooth functions.

The initial motivation for obtaining such a theorem came from trying to answer some questions
about generic smooth control systems; the potential usefullness of such result comes from the fact that o-minimality 
implies strong regularity properties of definable sets and maps ([vdD2]). In section 4, another possible application
is sketched.\n

To state our theorem precisely we introduce some notations. Fix $n\in\Z^+$, and let 
$\widehat{D^{\a}f}$, $\a=\a_1..\a_n\in (\Z^{\geq 0})^n$, be $n$-ary function symbols
indexed by the multiindex $\a$. Let ${\cal L}_{\widehat{Df}}$ denote the language obtained 
by adjoining the symbols $\widehat{D^{\a}f}$, $\a=\a_1..\a_n\in (\Z^{\geq 0})^n$, to the language ${\cal L}=\{0,1,+,-,\cdot,<\}$. 
Let $f\in C^{\infty}(\R^n,\R)$. We make $\R$ into ${\cal L}_{\widehat{Df}}$-structure
by giving the symbols $0,1,+,-,\cdot,<$ the usual
interpretation, and by interpreting $\widehat{D^{\a}f}$, $\a=\a_1..\a_n\in (\Z^{\geq 0})^n$, as functions which are equal to 
the corresponding partial derivatives of the given function $f$ on $[-1,1]^n$, and are equal to zero outside $[-1,1]^n$. We denote this structure and
its theory by ${\R}_{\widehat{Df}}$, $T_{\widehat{Df}}$ respectively. We call a subset of a Baire topological space $X$ residual if
it contains a countable intersection of open dense subsets of $X$. Our result is\n

{\bf Theorem A.} {\it There exists a residual subset $R\subset C^{\infty}(\R^n,\R)$, 
such that for each $f\in R$, the theory $T_{\widehat{Df}}$ is model complete and o-minimal.}\n

Moreover, let $\widehat{D^{\a}f_1},\widehat{D^{\a}f_2},..$ be $n_1,n_2,..$-ary respectively
function symbols for each  $\a=\a_1..\a_n\in (\Z^{\geq 0})^n$, and denote by ${\cal L}_{(\widehat{Df_i})_i}$,
the language obtained by adjoining to $\{0,1,+,-,\cdot,<\}$ these function symbols. 
Let $f_i\in C^{\infty}(\R^{n_i},\R)$, $i=1,2,..$ . We make $\R$ into
${\cal L}_{(\widehat{Df})_i}$-structure, and we denote by $T_{(\widehat{Df})_i}$ the theory of this structure.
One may generalize the proof of Theorem A and obtain\n

{\bf Theorem B.} {\it There exists a residual subset $R\subset \prod_i C^{\infty}(\R^n_i,\R)$, 
such that $\forall (f_i)_i\in R$, $T_{(\widehat{Df})_i}$ is model complete and o-minimal.}\n
 
Let us make a comment on the proof. In section 5 of [IY], a related problem was considered, 
roughly corresponding to estimating
the complexity of preimages of semialgebraic sets under jet extensions of generic smooth maps.
For their purposes, the authors construct a special equisingular Whitney stratification of such preimage,
and obtain the estimate sought from a certain numerical characteristic of this stratification. 
This characteristic is shown to be well defined by using Thom's First Isotopy Lemma.
We use here similar considerations, but develop them further: it is first established that for a generic function $f$, 
subsets which are defined by quantifier free ${\cal L}_{\widehat{Df}}$ formulas admit equisingular, in some sense,
Whitney stratifications. Then we use the Isotopy Lemma to show that projections of such sets admit
cylindrical decomposition.\n

We now outline the structure of the article. In section 1 we introduce Whitney stratifications 
of subsets of $\R^k$, whose strata project with constant rank on each {\it coordinate} subspace; 
we call such stratifications monotonic. In section 2
we study a certain subclass of quantifier free $\emptyset$-definable subsets from which each quantifier 
free $\emptyset$-definable set can be obtained by projection.  
For each set from this subclass we construct, using [GrY] and a construction of Thom-Boardman type, a monotonic Whitney
stratification whose strata are again sets from this subclass. This result allows us to construct a  
cylindrical decomposition for projections of such sets in section 3.  
As a corollary we obtain model completeness together with o-minimality, proving Theorem A. 
Finally, in section 4 we comment on the generalization to countably many generic smooth functions,
and sketch possible applications.\n

We have used above, and will be using below, some elementary logic-theoretic notions. For these, we refer the reader to
the very readable and concise introductory notes [Ch], whose terminology we follow here, or to any other 
introductory text on the topic. For more details on o-minimal structures see [vdD2]. Among articles which could be 
relevant to the theme discussed here we mention [KM], [W1], [W2], [RSW].\n\n

{\bf 1. Compact Whitney stratified subsets of $\R^n$ with monotonic strata.}\n

Below we summarize some facts about Whitney stratifications (for more details and proofs see 
[Ma1], [Ma2] or [GWPL], Chapters 1 and 2; the exposition here follows [GWPL]). 
The adjective 'smooth' will mean $C^{\infty}$, though we need only $C^1$ or sometimes $C^2$ smoothness
for the facts stated below. 
By a submanifold of $\R^n$ we mean an embedded smooth submanifold, not necessarily
connected, all of whose components have the same dimension. Let $X,Y\subset \R^n$ be submanifolds;
$X$ is said to be {\it Whitney regular} over $Y$ at $y\in Y$, 
if for any two sequences $(x_i)$ in $X$ and $(y_i)$ in $Y$, for which $x_i\neq y_i \ \forall i$, 
$x_i\rw y$, $y_i\rw y$, such that the sequences of lines $\overline{x_iy_i}$ and
of tangent spaces $T_{x_i}X$ converge, we have $\lim  \overline{x_iy_i} \subset \lim T_{x_i}X$.
$X$ is said to be Whitney regular over $Y$ if $\forall y\in Y$, $X$ is Whitney regular over $Y$ at $y$.
We may in fact take instead of $\R^n$ an arbitrary smooth manifold, and define Whitney regularity via
local charts.
If $X$ is Whitney regular over $Y$, then either $\overline{X-Y}\cap Y=\emptyset$ or $dim(Y)<dim(X)$.
A {\it stratification} ${\cal M}$ of $M\subset \R^n$ is a locally finite partition of $M$ into 
submanifolds, called strata.  It is said to be a Whitney stratification if each stratum is Whitney regular 
over any other stratum. We say that a stratified set has dimension $m$, if $m$ is the maximal 
dimension of its strata. We denote by ${\cal M}^q$ the subset of strata of $\cal M$ of dimension
smaller or equal to $q$, and denote the union of such strata by $M^q$.

The Cartesian product of two Whitney stratified subsets is again a Whitney stratified subset.
The stratification ${\cal M}$ is said to satisfy the {\it frontier condition} if $\forall X,Y\in{\cal M}, \ \overline{X}\cap Y\neq\emptyset$
implies $\overline{X}\supset Y$. 
A map is said to be transversal to ${\cal M}$ if
it is transversal to each stratum of ${\cal M}$. In this case, if ${\cal M}$ is Whitney, the preimages of 
the strata constitute again a Whitney stratification.

Any ($\emptyset$-definable) semialgebraic set admits a finite 
Whitney stratification with ($\emptyset$-definable) semialgebraic strata.\n

Our key analytical tool is Thom's Fisrt Isotopy Lemma. Let $M$ be a subset of a smooth manifold $N$
with a stratification ${\cal M}$. Let $f:N\rightarrow P$ be a smooth map into another smooth manifold $P$.
We say that $(M,{\cal M})$ is (topologically) trivial over $P$ if there exist a stratified set $F$ with a 
stratification ${\cal F}$, and
a homeomorphism $h:M\cong P\times F$, such that the following holds: 
each stratum of $M$ is sent to a stratum of $P\times{\cal F}$, 
and $\pi\circ h = f$, where $\pi:P\times F\rw P$ is the projection. We say that $(M,{\cal M})$ is (topologically) locally trivial  
over $P$ if each $p\in P$ has a neighbourhood $V$, such that $f^{-1}(V)\cap {\cal M}$ is again a stratification,
and $(f^{-1}(V)\cap M, f^{-1}(V)\cap {\cal M})$ is trivial over $f^{-1}(V)$.\n

{\bf Thom's First Isotopy Lemma.} ([GWPL], Chapter II, Theorem 5.2) {\it 
Let $(M,{\cal M})$ be a locally closed Whitney stratified subset of the smooth manifold $N$, and let $f:N\rw P$
be a smooth map such that for each $S\in {\cal M}$, $f|_S$ is a submersion and $f|_{\overline{S}\cap M}$ is a proper map.
Then $(M,{\cal M})$ is locally trivial over $P$.\n
} 

We will sometimes make use of the following fact:\n 

{\bf Proposition 1.1 } ([GWPL], Chapter II, Theorem 5.6, Corollary 5.7) 
{\it Let $N$ be a smooth manifold, and let $M\subset N$ be a closed set with a Whitney stratification $\cal M$.
Then the components of strata of $\cal M$ form another Whitney stratification, which moreover
satisfies the frontier condition.}\n

Let the coordinates on $\R^n$ be denoted by $x_1,..,x_n$.
By a {\it coordinate plane} in $\R^n$ we mean 
any set of the form $\{(x_1,..,x_n)\in\R^n :x_{i_1}=0,..,x_{i_m}=0\}$,
for some choice of $m\leq n$  and $1\leq i_1<..<i_m\leq n$.  
We say that a submanifold $S\subset\R^n$ is {\it monotonic}, if the
projection of $S$ on any coordinate plane in $\R^n$ is a map of constant rank.
We call a stratification of a subset of $\R^n$ monotonic if all
its strata are monotonic. Let $M\subset\R^k$ be a set with a monotonic stratification $\cal M$,
and let $P$ be a coordinate plane in $\R^k$. We denote by $M(q,P)$ the union of strata of $M$ which project on $P$ with
rank smaller or equal to $q$, and by ${\cal M}(q,P)$ the corresponding stratification.\n

{\bf Proposition 1.2.} {\it
Let $M\subset\R^k$ be a closed Whitney stratified set and let $q\geq 0$ be an integer. Then $M^q$ 
is a closed set. Suppose further that
the strata of $\cal M$ are monotonic, and let $P$ be a coordinate plane in $\R^k$. 
Then the set $M(q,P)$ is closed as well.
}

{\bf Proof.} The set $M^q$ is closed since $M$ is closed and the closure of a stratum of dimension $d$, in a Whitney
stratification, can only intersect strata of dimension smaller than $d$. The set $M(q,P)$ is closed since $M$ is closed
and because of the following property of Whitney stratifications. Namely, suppose $y_i$ is a sequence of points on a stratum
$S$ converging to a point $x$ on a stratum $S'$, and suppose that the sequence $T_{y_i}S$ converges to a subspace $T$. 
Then $T_xS'\subset T$.
\ \ \ \ $\Box$\n 

{\bf Lemma 1.3.} {\it Let $S$ be a monotonic submanifold of $\R^n$, let $\pi:\R^n\rw\R^k$ denote the projection
on a coordinate plane $\R^k\subset\R^n$, and let $x\in\R^k$. Then the set $\pi^{-1}(x)\cap S$ is a monotonic
submanifold of $\R^n$.}

{\bf Proof.} Basically, an exercise in linear algebra.
\ \ \ $\Box$\n

{\bf Lemma 1.4.} {\it Let $M\subset\R^n$ be a compact set with a monotonic Whitney stratification $\cal M$ whose strata are connected. 
Then the closure of each stratum contains a $0$-dimensional stratum.}

{\bf Proof.} Let $S\in{\cal M}$, $dim(S)>0$. The stratum $S$ cannot be closed, since then it would be compact, and its projection on
any coordinate axis would have a critical point. This would imply, since $S$ is monotonic, that the rank of each such projection is zero,
and consequently $dim(S)=0$, which contradicts our assumption. 
The set $M$ is closed, hence $\ov{S}$ must intersect another stratum $S'$ of $\cal M$. 
Since $\cal M$ is Whitney, $dim(S')<dim(S)$. 
Since the strata of $\cal M$ are connected, $\cal M$ satisfies the frontier condition (Proposition 1.1), and thus $S'\subset\ov{S}$.
The submanifold $S'$ is monotonic and $\ov{S'}\subset\ov{S}$, so we may repeat the argument and eventually conclude 
that $\ov{S}$ contains a $0$-dimensional stratum. \ \ \ \ $\Box$\n   

The following two lemmas about monotonic Whitney stratifications will be used in section 3.\n

{\bf Lemma 1.5.} {\it
Let $M\subset\R^n$ be a compact set with a monotonic Whitney stratification $\cal M$ whose strata are connected, 
and let $\pi:\R^n\rw \R$ be the projection on a
one dimensional coordinate plane of $\R^n$. Let ${\cal M'}$ be a subset of $\cal M$, and let $M'$ be the union of 
strata of ${\cal M}'$. Then there exists a set ${\cal T}$ of $1$-dimensional strata of ${\cal M}$,
and a set $\cal P$ of $0$-dimensional strata of ${\cal M}$, such that the following is true:\n
i) for each $S\in{\cal T}$, $rank(\pi_S)=1$,\n
ii) the boundary of $\pi(S)$, for each each $S\in\cal T$, is contained in $\pi(M^0)$,\n
iii) the projections of $\cup_{S\in{\cal T}}S$ and $\cup_{S\in{\cal P}}S$ are disjoint and their union is equal to $\pi(M')$.  
}

{\bf Proof.} Let $C=\pi(M^0)$. For each $p\in\R-C$, $\pi^{-1}(p)\cap M$
has the monotonic Whitney stratification ${\cal M}_p=\{\pi^{-1}(p)\cap S:S\in{\cal M}\}$,
whose strata are monotonic by Lemma 1.3. By Lemma 1.4, $p\in \pi({\cal M}_p^0)$.
Therefore for each $S'\in{\cal M}'$, $\pi(S')-C$ is contained in $\cup\pi(L)$ where $L$ ranges over
$1$-dimensional strata of $\ov{S'}$. Since the strata of $\cal M$ are connected, $\cal M$ satisfies the frontier condition (Proposition 1.1),
and thus $\pi(S')-C$ is in fact equal to  $\cup\pi(L)-C$. If we take ${\cal T}$ to consist of $L\in {\cal M}$ such that
$dim(L)=1$, $L\subset\ov{M'}$, $rank(\pi|_L)=1$, and we take ${\cal P}$ consist of $P\in{\cal M}^0$ such that
$\pi(P)\in \pi(M')-\cup_{L\in{\cal T}}\pi(L)$, we see that i) and iii) hold. Using Lemma 1.4, one shows that 
ii) holds as well. 
\ \ \ \ $\Box$\n

Let $f:N\rw P$ be a smooth map, and let $M\subset N$ have a Whitney stratification $\cal M$. We say that $p\in P$
is a regular value of $f|_M$, if $p$ is a regular value of each $f|_S$, $S\in{\cal M}$.\n

{\bf Lemma 1.6.} {\it  Let $\pi:\R^n\rw \R^k$ be the 
projection of a compact set $M\subset\R^n$ to a coordinate plane  of $\R^n$. 
Suppose that $M$ admits a monotonic Whitney stratification.
Then the number of components of $\pi^{-1}(x)\cap M$ is
uniformly bounded over $x\in\R^k$.}

{\bf Proof}  
Choose some monotonic Whitney stratification of $M$ and denote it by $\cal M$.
By taking, if necessary, the connected components of strata, 
we may assume that the strata are connected and thus $\cal M$ satisfies the frontier condition (Proposition 1.1). 
Let $A\subset M$. We denote by $N_x(A)$, $x\in\R^k$ the number of components of 
$\pi^{-1}(x)\cap A$. Observe that
$$
N_x(M)\leq \sum_{S\in\M} N_x(\ov{S}).
$$  
Take a stratum $S$ and denote the rank with which it projects to $\R^k$
by $r(S)$.
By Lemma 1.3, $\pi^{-1}(x)\cap S$ is a monotonic submanifold of $\R^n$,
and thus by Lemma 1.4 none of its components can be compact unless $dim(S)=r(S)$.
If $dim(S)>r(S)$, each component of $\pi^{-1}(x)\cap \ov{S}$ must therefore 
intersect the {\it frontier} of $S$, defined as $\fr(S)=\ov{S}-S$. 
We conclude that in the case $dim(S)>r$,
the number of components of $\pi^{-1}(x)\cap \ov{S}$
is bounded by the number of components of 
$\pi^{-1}(x)\cap \fr(S)$, which, by the frontier condition, is itself
bounded by $\sum_{S'\in \fr(S)}N(\ov{S'})$. The dimension of 
$S'\in \fr(S)$
must be strictly smaller than $dim(S)$ (since the stratification is 
Whitney). Repeating the argument as many times as needed, 
we may bound $N_x(M)$ by a sum of $N_x(\ov{S})$ over strata $S$ 
for which $dim(S)=r(S)$ (each term may appear more than once in this sum). 

Thus, to prove that the number of components
of $\pi^{-1}(x)\cap M$ is uniformly bounded over $x\in\R^k$, 
it is sufficient to show, that for each $S\in\M$ for which $dim(S)=r(S)$, the cardinality
of $\pi^{-1}(x)\cap S$ is uniformly bounded over $x$. Since one may find a
coordinate plane in $\R^k$ of dimension $dim(S)$, on which $S$ projects with
rank $dim(S)$, we may assume that $S$ maps to $\R^k$ by a local diffeomorphism.

Let $\R^k=\R^{k-1}\times\R$, and consider the projection $\pi':\R^n\rw\R^{k-1}$.
Denote by $E\subset\R^{k-1}$ the image
of strata of $\ov{S}$ which map with rank smaller than $k-1$ to $\R^{k-1}$, 
and denote by $\pi'':\R^k\rw \R^{k-1}$ the projection from $\R^k$ to $\R^{k-1}$. 
For all $y\in\R^{k-1}-E$, $(\pi')^{-1}(y)\cap \ov{S}$ is then a 
one dimensional Whitney stratified set, with monotonic strata, which we denote by $(\ov{S})_y$. 
It looks like a graph, with vertices and edges (with possibly multiple edges connecting 
a pair of vertices), each edge being a one dimensional submanifold. 
The number of vertices of the graph is bounded by 
the number of regular preimages of $y$ in $(\ov{S})^{k-1}$ under $\pi'$.
It is a consequence of Thom's First Isotopy Lemma, that the maximal degree $D$ of the vertices 
is bounded uniformly over all $y\in\R^{k-1}-E$ (see the proof of Theorem 5.1 in [IY]).  
Moreover, each edge of the graph
must land in at least one vertex of the graph, since the edges are monotonic and bounded.
Let $L$ be an edge which
projects with rank 1 to the line $(\pi'')^{-1}(y)$.
There can be only one preimage of $x\in (\pi'')^{-1}(y)$ in $L$, since otherwise there
would be a point on L with tangent which is parallel to the plane $\pi^{-1}(x)\subset(\pi')^{-1}(y)$, 
contradicting monotonicity of $L$. 

Denote by $Z_y\subset (\pi'')^{-1}(y)$ the finite set
of projections of vertices of the graph. 
Thus, for all $x\in(\pi'')^{-1}(y)-Z_y$, $y\in\R^{k-1}-E$, the number of preimages of $x$ in $\ov{S}$ 
is bounded by the product of
the maximal graph degree $D$ and the number of regular preimages of $y\in\R^{k-1}$ in $(\ov{S})^{k-1}$.
Suppose that the latter is uniformly bounded over $\R^{k-1}-E$. Then there exists a dense set $R\subset\R^k$, 
such that the number of preimages of $x$ in $\ov{S}$ is 
uniformly bounded over $x\in R$, say by $N>0$. If there existed $x_0\in\R^k$ with more than $N$ regular preimages,
this would be true for every $x$ in a sufficiently small neighbourhood of $x_0$, thus contradicting the density
of $R$. Thus, if the number of regular preimages of $y\in\R^{k-1}$ in $(\ov{S})^{k-1}$ is uniformly bounded, 
then also the number of regular preimages of $x\in\R^{k}$ in $\ov{S}=(\ov{S})^k$ 
is uniformly bounded. Proceeding by induction, we conclude that the number of 
regular preimages of $x\in\R^{k}$ in $\ov{S}=(\ov{S})^k$ 
is indeed uniformly bounded over $\R^k$. 

By the reduction made, this proves the lemma.\ \ \ \ \ $\Box$\n\n

{\bf 2. Basic $f$-sets and their monotonic Whitney stratifications.}\n

In this section we study the quantifier free $\emptyset$-definable sets of the structure ${\R}_{\widehat{Df}}$. 
The class of sets which we now introduce, the $\emptyset$-definable {\it basic $f$-sets}, will be shown to form a 
subclass of the class of $\emptyset$-definable quantifier free sets, with the property
that each $\emptyset$-definable quantifier free set is a projection of a \bef. 
It will be also shown that each $\emptyset$-definable basic $f$-set
has a finite Whitney stratification, whose strata are monotonic and are $\emptyset$-definable basic $f$-sets themselves.
We use this in section 3 to establish a cylindrical decomposition result for projections
of $\emptyset$-definable quantifier free sets.

Below, $J^m(\R^n,\R)$ denotes the space of $m$-jets of functions from $\R^n$ to $\R$, and
$j^m f(x)$ denotes the $m$-th jet of a smooth function $f:\R^n\rw\R$ at $x\in\R^n$.
We denote by $J(r,m,n)$ the $r-th$ power of $J^m(\R^n,\R)$, which we may 
identify with the space of all tuples $(j^m f(x_1),..,j^m f(x_r))$, $f\in C^{\infty}(\R^n,\R)$, $(x_1,..,x_r)\in(\R^n)^r$.
We call it the $(r,m)$-{\it multijet} space (compare with [GG], page 57; 
our definition is different since we {\it do not} require the diagonal in
$(\R^n)^r$ to be excluded). We call the map $j^{r,m}f:(\R^n)^r\rw J(r,m,n)$, defined by
$(x_1,..,x_r)\mapsto (j^m f(x_1),..,j^m f(x_r))$, the multijet extension
of $f$. For our purposes we identify $J(r,m,n)$ with a Euclidean space 
of corresponding dimension. 

One encounters multijet preimages of $\emptyset$-definable semialgebraic sets 
among the quantifier free $\emptyset$-definable sets of the structure ${\R}_{\widehat{Df}}$ (for example, 
let $f\in C^\infty(\R,\R)$, and consider the set $\{(x_1,x_2)\in [-1,1]^2:f(x_1)>f(x_2)\}$). 
Since our aim is to construct special Whitney stratifications for quantifier free $\emptyset$-definable sets,
we at least should be able to produce Whitney stratifications for generic multijet preimages
of semialgebraic sets. If we try to use the Multijet Transversality
Theorem as stated in [GG] (Chapter II, Theorem 4.13), 
we see that it only implies that for a generic function $f$ the multijet preimage of
a semialgebraic set in $(\R^n)^r-Diag$ is Whitney stratifiable, 
and not that the preimage in  $(\R^n)^r$ itself is Whitney stratifiable (we denote by $Diag$ the diagonal). 
This particular problem was addressed in [GrY], the result of which we state below.
We need to introduce first the notion of divided differences (for more details see, for instance, [BZ]).\n

{\bf Definition 2.1.} Let $f\in C^\infty(\R,\R)$, and let $k\in \Z^{\geq 0}$. The {\it divided difference of order $k$ of $f$}, which we
denote by $\Delta^k f$, is a function from $\R^{k+1}$ to $\R$ defined as follows. For $k=0$, $\D^0f=f$. 
For $k>0$, let $(x_1,..,x_{k+1})\in\R^{k+1}$ be such that $x_i\neq x_j$ $\forall i\neq j$.
Then $\D^kf(x_1,..,x_{k+1})$ is defined by the following recursion relation:
$$
\D^kf(x_1,..,x_{k+1})=(\D^{k-1}f(x_1,..,x_{k})-\D^{k-1}f(x_2,..,x_{k+1}))/(x_1-x_{k+1}).
$$

It can be shown that this defines, by continuity, a smooth function from $\R^{k+1}$ to $\R$, which is 
in fact symmetric. There exists an explicit formula, easily proved by induction, 
to compute the divided difference at points lying on the 
diagonal, i.e. at points $(x_1,..,x_{k+1})$ for which there exist $i,j$, $i\neq j$, such that $x_i= x_j$.
For stating this formula, we introduce the following notation. Let $m$ be a nonnegative integer; 
denote by $\d^m:\R\rw\R^{m+1}$ the map which sends $x\in\R$ to the tuple $(x,..,x)\in\R^{m+1}$, 
all of whose entries are equal to $x$.\n

{\bf Proposition 2.2.} Let $f\in C^\infty(\R,\R)$. Then 
$$
\D^{m_1+..+m_l+l-1}f(\d^{m_1}(u_1),..,\d^{m_l}(u_l))={1\over m_1!..m_l!}D^{m_1..m_l}(\D^{l-1}f(u_1,..,u_l)).
$$

{\bf Remark.} Note that since divided differences are symmetric functions of their arguments, one may indeed use 
Proposition 2.2 to compute $\D^kf$ at any given point of $\R^{k+1}$. Moreover, Proposition 2.2 implies that
there exist $t>k$ and a semialgebraic function $\rho:J(k,t,1)\rw \R$, so that $\D^kf=\rho\circ j^{k,t}f$.\n

One may introduce divided differences also for functions of several variables. 
Let $f\in C^{\infty}(\R^n,\R)$, $f=f(x_1,..,x_n)$, and fix $1\leq i \leq n$. 
We denote by $\D_{x_i}^kf$ the function which we obtain by 
taking the divided difference of $f$ of order $k$ w.r.t. the variable $x_i$, while keeping the other variables fixed.
Denote the variables on which $\D_{x_i}^kf$ depends by $x_1,..,x_{i-1},x_{i1},..,x_{i \ k+1},x_{i+1},..,x_n$.
We may now take a divided difference w.r.t. another variable $x_j$, $j\neq i$,
and repeat this operation w.r.t. other variables.\n 

{\bf Definition 2.3.} Let $f\in C^\infty(\R^n,\R)$, $f=f(x_1,..,x_n)$, and let $(\a_1,..,\a_n)\in (\Z^{\geq 0})^n$. 
The divided difference of order $(\a_1,..,\a_n)$ of $f$, denoted 
by $\Delta^{\a_1..\a_n} f$, is the function defined as 
$$
\Delta^{\a_1..\a_n} f=\D^{\a_1}_{x_1}\circ ..\circ \D^{\a_n}_{x_n}f.
$$

The function $\Delta^{\a_1..\a_n} f$ depends on $\a_1+..+\a_n+n$ variables, which we denote by  $x_{11}, .., x_{1 \a_1+1},
$ $... \ ,$ $x_{n1}, .., x_{n \a_n+1}$. We make the following convention: $\Delta^{\a_1..\a_n} f (X_1,..,X_n)$, where $X_1,..,X_n$ are
tuples of real numbers, means that we substitute for $x_{i1}, .., x_{i \a_i+1}$ the first $\a_i+1$ entries of the tuple $X_i$,
for each $1\leq i\leq n$.

Let $m,r$ be nonnegative integers, $r>0$. Denote by $diag^m:\R^r\rw \R^{(m+1)r}$ the map which sends $(x_1,..,x_r)$
to $(\d^m(x_1),..,\d^m(x_r))$, and by $\widetilde{\cal D}^{r,m}f:(\R^n)^r\rw \R^{(r(m+1))^n}$ the map which sends
$(x_{11},..,x_{1r},...,x_{n1},..,x_{nr})$ to the collection of 
$$
\Delta^{\a_1..\a_n} f(diag^m(x_{11},..,x_{1r}),..,diag^m(x_{n1},..,x_{nr})),
$$
ordered lexicographically on multiindices, where $0\leq \a_i\leq r(m+1)-1$ for each $1\leq i \leq n$. Denote 
by ${\cal D}^{r,m}f: (\R^n)^r\rw \R^{nr+(r(m+1))^n}$ the map which sends $u\in(\R^n)^r$ to $(u,\widetilde{\cal D}^{r,m}f(u))$, and call it
the divided difference extension of $f$. Below, we denote by $\wt{D}(r,m,n)=\R^{(r(m+1))^n}$ the target space of $\wt{{\cal D}}^{r,m}f$,
and by ${D}(r,m,n)=\R^{nr+(r(m+1))^n}$ the target space of ${{\cal D}}^{r,m}f$.\n

{\bf Theorem 2.4 [GrY].} {\it \n
i) Let $M\subset D(r,m,n)$ be a submanifold.
There exists a residual set $R\subset C^{\infty}(\R^n,\R)$, such that $\forall f\in R$, ${\cal D}^{r,m}f$ is transversal to $M$.\n
ii) there exists a polynomial map with rational coefficients $\b:D(r,m,n)\rw J(r,m,n)$, such that $j^{r,m}f=\pi\circ {\cal D}^{r,m}f$.\n
} 

Let us write down a simple example. Take $f\in C^{\infty}(\R,\R)$; then $(x_{11},x_{12})$ is mapped by $j^{2,0}f$ to $(x_{11},f(x_{11}),x_{12},f(x_{12}))$, and 
by ${\cal D}^{2,0}f$ to $(x_{11},x_{12},f(x_{11}),(f(x_{11})-f(x_{12}))/(x_{11}-x_{12}))$. In this case the map $\b$ is a birational morphism; this 
is not true for $n\geq 2$.

Note that the theorem implies that for any semialgebraic set $S\subset D(r,m,n)$ there exists a residual subset $R\subset C^{\infty}(\R^n,\R)$,
such that for each $f\in R$, $(j^{r,m}f)^{-1}(S)$ has a Whitney stratification. (Indeed, choose a Whitney stratification of
$\b^{-1}(S)$. According to the theorem there exists a 
residual set $R\subset C^{\infty}(\R^n,\R)$, such that for each $f\in R$, ${\cal D}^{r,m}f$ is transversal to each of the strata of 
this Whitney stratification.
Now recall that transversal preimages of Whitney stratified sets have a natural Whitney stratification.) 
Just having a Whitney stratification, however,
is not sufficient, as we would also like the strata to be monotonic. It is this requirement which does not allow us  (unless $n=1$) to limit ourselves
to the subclass of multijet preimages of semialgebraic sets, if we wish the strata to be again sets from this subclass.
It turns out that a suitable way to proceed is to consider instead the subclass 
of preimages of semialgebraic sets under divided difference extensions of $f$. More precisely:\n

{\bf Definition 2.5. } Let $f\in C^{\infty}(\R^n,\R)$.
We call $M\subset\R^p\times (\R^n)^r$, $p\geq 0$, $n,r\geq 1$, a $\emptyset$-definable basic $f$-set, if there exists $m\geq 0$, and a
$\emptyset$-definable semialgebraic set $A\subset \R^p\times D(r,m,n)=\R^p\times(\R^n)^r\times \wt{D}(r,m,n)$, such that\n
i) $M=(id\times {\cal D}^{r,m}f)^{-1}(A)$,\n 
ii) $A$ is a subset of $\R^p\times[-1,1]^{nr}\times \wt{D}(r,m,n)$.\n

Here, the map $id\times {\cal D}^{r,m}f$ maps $(u,v)\in \R^p\times (\R^n)^r$ to $(u,{\cal D}^{r,m}f(v))$.
The requirement ii) corresponds to restricting $f$ to the unit cube; it implies that $M$ is a subset of $\R^p\times[-1,1]^{nr}$.\n

{\bf Remark.} Note that the numbers $p,n,r$ are not determined by the dimension of the ambient space of $M$.
When we speak below of \bef s inside some ambient space $\R^q$, we always provide a representation of $\R^q$ in the form 
$\R^p\times(\R^n)^r$, from which these numbers can be determined.\n

Let $\widetilde{w}^{r,m}f:(\R^n)^r\rw \widetilde{W}(r,m,n)$, where $\widetilde{W}(r,m,n)$ is a Euclidean space, 
denote the map whose components are $D^{\a_1..\a_n}f(x_{1i_1},..,x_{ni_n})$, $ \ 0\leq \a_1,..,a_n\leq t$, $1\leq i_1,..,i_n\leq r$.
Denote by $w^{r,t}f$ the map which sends $u\in(\R^n)^r$ to $(u,\widetilde{w}^{r,t}f(u))\in(\R^n)^r\times \widetilde{W}(r,t,n)$. Denote 
$(\R^n)^r\times \widetilde{W}(r,t,n)$ by $W(r,t,n)$.\n

{\bf Lemma 2.6.} {\it Let $f\in C^{\infty}(\R^n,\R)$, and let $m\geq 0, r>0$.\n 
i) There exist $t\geq m$ and a $\emptyset$-definable semialgebraic map\n
$\rho_{r,t,n}:W(r,t,n)\rw D(r,m,n)$,
such that ${\cal D}^{r,m}f = \rho_{r,t,n}\circ w^{r,t}f$.\n
ii) there exists a polynomial map with rational coefficients $\sigma: D(r,m,n)\rw W(r,m,n)$, so that
$w^{r,m}f=\sigma\circ {\cal D}^{r,m}f$.\n}  

{\bf Remark.} Note that although the definition of ${w}^{r,m}f$ may seem similar to that of
$j^{r,m}f$, it is not the same map, and in fact
part ii) of Lemma 2.6 is stronger than part ii) of Theorem 2.4. 
The semialgebraic map $\rho_{r,t,n}$ whose existence is asserted in part i), is not 
continuous on all of $W(r,t,n)$.\n

{\bf Sketch of the proof.} Part i) is implied by Proposition 2.2 and by symmetricity of divided differences of functions of
one variable. Part ii) follows from the method of proof of Theorem 2.4 ii) in [GrY] (see [GrY], Remark on pg. 359). 
\ \ \ \ $\Box$\n

{\bf Lemma 2.7.} {\it Let $f\in C^{\infty}(\R^n,\R)$. 
Then each $\emptyset$-definable basic $f$-set is a quantifier free $\emptyset$-definable set of
$\R_{\widehat{Df}}$.
}

{\bf Sketch of the proof.}  
Take a $\emptyset$-definable basic $f$-set $M$. By definition, there exist nonnegative integers $r,m,n,p$
and  a $\emptyset$-definable semialgebraic set $A\subset \R^p\times D(r,m,n)$, such that $M=(id\times {\cal D}^{r,m}f)^{-1}(A)$, and such that
$A$ is a subset of $\R^p\times[-1,1]^{nr}\times\wt{D}(r,m,n)$. 
By Lemma 2.6 i) there exist $t\geq 0$ and a $\emptyset$-definable semialgebraic map $\rho_{r,t,n}:W(r,t,n)\rw D(r,m,n)$, such that
${\cal D}^{r,m}f = \rho_{r,t,n}\circ w^{r,t}f$. 
Thus it is also true that $M=((id\times \rho_{r,t,n})\circ(id\times w^{r,t}f))^{-1}(A)$, hence
$M=(id\times w^{r,t}f)^{-1}( (id\times \rho_{r,t,n})^{-1}(A))$.

If we show that $C=(id\times \rho_{r,t,n})^{-1}(A)$ is a subset of $\R^p\times [-1,1]^{nr}\times\wt{W}(r,t,n)$, 
we are done, since this allows us to write 
a formula in ${\cal L}_{\widehat{Df}}$ for the set $M$ (indeed, $C$ is a $\emptyset$-definable semialgebraic set, and
$ (id\times w^{r,t}f)^{-1}(C)$ can then be defined by a quantifier free formula involving $\widehat{D^{\a}f}$). But since  
${\cal D}^{r,m}f = \rho_{r,t,n}\circ w^{r,t}f$, $\rho_{r,t,n}$ is the identity map on the first $nr$ coordinates. Hence, 
the fact that $A$ is a subset of $\R^p\times[-1,1]^{nr}\times \wt{D}(r,m,n)$ implies that 
$(\rho_{r,t,n})^{-1}(A)\subset \R^p\times [-1,1]^{nr}\times \wt{W}(r,t,n)$. \ \ \ \ $\Box$\n

We wish to show that all quantifier free $\emptyset$-definable sets of $\R_{\widehat{Df}}$ are in fact 
projections of $\emptyset$-definable basic $f$-sets (and therefore the class of 
projections of quantifier free $\emptyset$-definable sets is identical to the class of projections of  $\emptyset$-definable basic $f$-sets).
We will need the notion of {\it depth} of quantifier free formulas, which we now define for 
an arbitrary given language $L$. 
We say that a term has depth $0$ if it is a variable or a constant of $L$.
Consider the term $f(t_1,..,t_s)$, where $f$ is a s-ary function symbol of $L$, and $t_1,..,t_s$ are
also terms. We say that $f(t_1,..,t_s)$ has depth $i$, if
the maximal depth of the terms $t_1,..,t_n$ is equal to $i-1$. The depth of an atomic formula is defined as the maximal depth
of the terms on which it depends. We say that a quantifier free formula, which is a boolean combination of atomic formulas,
has depth $i$, if $i$ is the maximal depth of the atomic formulas of which the boolean combination is formed.\n

{\bf Lemma 2.8.} {\it Let $f\in C^\infty(\R^n,\R)$. Then each quantifier free $\emptyset$-definable set of $\R_{\widehat{Df}}$
is the projection of a $\emptyset$-definable basic $f$-set on a coordinate plane.}

{\bf Sketch of the proof.} One shows using induction, that for each quantifier free formula
$\phi(u_1,..,u_s)$ of any language $L$ there exists an equivalent existential formula
$$
\exists v_1 .. \exists v_t ( \psi(u_1,..,u_q,v_1,..,v_t),
$$
where $\psi(u_1,..,u_q,v_1,..,v_t)$ is a quantifier free formula  of depth at most $1$ (i.e. not involving
compositions of functions).
Moreover, $\psi$ can be so chosen that all its terms of depth $1$ do not have constants as
arguments, and all depend on pairwise disjoint groups of variables. Applying this observation to
the language ${\cal L}_{\widehat{Df}}$, one further shows that in fact the formula $\psi$ can be chosen in such a way, that it defines
the preimage of a $\emptyset$-definable semialgebraic set $S\subset \R^p\times J(r,m,n)$ under $id\times j^{r,m}f$, for some 
$r,m,n,p$. Moreover, if we denote by $\pi: J(r,m,n)\rw (\R^n)^r$ the projection on a subspace of $J(r,m,n)$ with the property that
$\pi\circ j^{r,m}f(u)=u \ \forall u\in(\R^n)^r$, then $\pi(S)\subset[-1,1]^{nr}$.  
By Theorem 2.4 ii), there exists a polynomial map $\b$ with rational coefficients from $D(r,m,n)$ to $J(r,m,n)$, such that
$j^{r,m}f=\b\circ {\cal D}^{r,m}f$. Thus 
$$
(id\times j^{r,m}f)^{-1}(S) \ = \ (id\times {\cal D}^{r,m}f)^{-1}( (id\times\b)^{-1}(S)).
$$
Since $\pi(S)\subset[-1,1]^{nr}$, $(id\times\b)^{-1}(S)$ is a subset of $\R^p\times[-1,1]^{nr}\times \wt{D}(r,m,n)$. 
Moreover, since $\b$ is a polynomial map with rational coefficients, $(id\times\b)^{-1}(S)$
is a $\emptyset$-definable semialgebraic set. Therefore, the quantifier free formula $\psi$ defines a 
$\emptyset$-definable basic $f$-set.\ \ \ \ $\Box$\n

We illustrate the steps which were not detailed in the proof of Lemma 2.8 in a simple case. 
Let $n=1$, and consider the set defined by the depth $2$ formula $\wf(\wf(x_1))>0$
(we write $\wf$ instead of the ${\cal L}_{\wh{Df}}$ symbol $\wh{D^0f}$). 
Note that this set can be defined by the following depth $1$ formula 
$$
\exists x_2 ( \wf(x_2)>0 \wedge \wf(x_1)=x_2),
$$
By further transformation, we arrive at the equivalent formula
$$
\exists x_2 \exists x_3 \exists x_4 ( \ \ \wf(x_3)>0 \wedge |x_3|\leq 1 \wedge |x_4|\leq 1 \wedge x_2=x_3 \wedge x_2=0\wedge |x_1|>1 \ \ \ \ \vee
$$
$$
\wf(x_3)>0 \wedge |x_3|\leq 1 \wedge |x_4|\leq 1 \wedge x_2=x_3 \wedge x_2=\wf(x_4) \wedge x_1=x_4 \ \ ),
$$
which is the preimage under $id\times j^{2,0}f$ of a semialgebraic subset of $\R^2\times [-1,1]^2\times\R^2$.\n

We now aim to show that for a generic smooth function $f\in C^\infty(\R^n,\R)$, 
every $\emptyset$-definable basic $f$-set has a Whitney monotonic stratification whose strata are 
$\emptyset$-definable basic $f$-sets. By 'generic' we mean that $f$ is such that
for all $p,r,m$, $id\times{\cal D}^{r,m}f$ is transversal
to any $\emptyset$-definable semialgebraic submanifold of $\R^p\times D(r,m,n)$. 
We wish to show that the set of such functions, which we denote by $R_n$, is residual.\n

{\bf Proposition 2.9.} {\it
For each $n\geq 1$, $R_n$ is a residual set.
}

{\bf Proof.} Fix the numbers $p,r,m$, and let $N\subset\R^p\times D(r,m,n)$ be a submanifold. It is a simple generalization of Theorem 2.4, that
the set $Q(p,r,m,N)$ of functions $f\in C^\infty(\R^n,R)$ for which $id\times{\cal D}^{r,m}f$ is transversal to $N$, is residual.
The set $R_n$ is the intersection of the sets  $Q(p,r,m,N)$, as $p,r,m$ range over nonnegative integers, and $N$ ranges over
$\emptyset$-definable semialgebraic submanifolds of $\R^p\times D(r,m,n)$. Since the family of all $\emptyset$-definable semialgebraic
submanifolds is countable, we conclude that $R_n$ is residual. \ \ \ $\Box$\n

We will also need the following two facts.\n

{\bf Proposition 2.10.} {\it
Let the set $M\subset\R^k$ have a Whitney stratification $\cal M$ of dimension $m$. 
Let $\cal A$ be a Whitney stratification of $M^q$, $q\leq m$, which refines ${\cal M}^q$. Then 
$({\cal M}-{\cal M}^q)\cup {\cal A}$, a stratification of $M$ which refines $\cal M$, is Whitney.}  

{\bf Proof.} If $X,Y$ are submanifolds and $X$ is Whitney over $Y$,
then $X$ is necessarily Whitney over any submanifold of $Y$. If moreover $X,Y$ are disjoint and $dim(X)\leq dim(Y)$ then $\ov{X}\cap Y=\emptyset$,
and therefore any submanifold of $X$ is Whitney over $Y$. The proof follows from these facts. \ \ \ $\Box$\n

{\bf Proposition 2.11.} {\it
Let $A_1,..,A_l\subset\R^k$ be $\emptyset$-definable semialgebraic sets. Then there exists 
a finite Whitney stratification ${\cal A}$ of $A_1\cup..\cup A_l$, whose strata are 
$\emptyset$-definable semialgebraic sets, such that for each stratum $T\in{\cal A}$ and $1\leq i\leq l$, either $T\subset  A_i$ or $T\cap A_i=\emptyset$.    
}

{\bf Proof.} This is a well known fact. \ \ \ \ $\Box$\n

{\bf Lemma 2.12.} {\it Let $f\in C^\infty
(\R^n,\R)$. Fix the numbers $r,m$, and consider the map
${\cal D}^{r,m}f:(\R^n)^r \rw D(r,m,n)$.
There exists $t>m$, such that each partial derivative of first order of any of the components of ${\cal D}^{r,m}f$
is a composition of a $\emptyset$-definable semialgebraic function
on $D(r,t,n)$ with ${\cal D}^{r,t}f:(\R^n)^r\rw D(r,t,n)$.
}

{\bf Sketch of the proof.} Let
$$
\Delta^{\a_1..\a_n} f(diag^m(x_{11},..,x_{1r}),..,diag^m(x_{n1},..,x_{nr})),
$$
($0\leq \a_j\leq r(m+1)-1$, $j=1,..,n$), be a given component of the map ${\cal D}^{r,m}f$. 
Apply Proposition 2.2 separately to each variable of $f\in C^\infty(\R^n,\R)$. This
allows us to represent the given component in the form:
$$
c \ \cdot \ D^{m..m l_n}_{x_{n1}..x_{n w_n+1}} \ \D^{w_n}_{x_n} \  \ ... \  \ D^{m..m l_i}_{x_{11}..x_{1 w_1+1}} \ \D^{w_1}_{x_1}f
$$ 
evaluated on $(x_{11},..,x_{1 \ w_1+1}, \ ... \ ,x_{n1},..,x_{n \ w_n+1})$. Here $(m+1)w_i+l_i=\a_i$, $0\leq l_i\leq m$, 
$i=1,..,n$, and $c$ is some rational constant.
Apply the partial derivative operator $\partial_{x_{pq}}$, $1\leq p\leq n$, $1\leq q\leq w_p+1$; it commutes with 
other partial derivative operators and with $\D^{w_n}_{x_n},..,\D^{w_{p+1}}_{x_{p+1}}$, so one may write the result
as
$$
c \ \cdot \ D^{m..m l_n}_{x_{n1}..x_{n w_n+1}} \ \D^{w_n}_{x_n} \ ... \  D^{m..(m+1)..m l_i}_{x_{p1}..x_{pq}..x_{1 w_1+1}} \ \D^{w_p}_{x_p}
    \ ...\ D^{m..m l_i}_{x_{11}..x_{1 w_1+1}} \ \D^{w_1}_{x_1}f.
$$
Note that this expression will {\it not} in general be in the form of a component of $\D^{r,t}f$ for some $t>m$. However, up to a multiplication
by a rational constant, it is a divided difference of $f$ (we apply again Proposition 2.2, this time in the opposite direction). 
Consequently, it can be shown (by the same argument which proves Lemma 2.6 i)) 
to be the composition of a $\emptyset$-definable semialgebraic function $\rho$, defined on $W(r,t,n)$, and of $w^{r,t}f$, for some $t>m$.  

By Lemma 2.6 ii), there exists a polynomial map with rational coefficents $\s:D(r,t,n)\rw W(r,t,n)$, so that
$\s\circ{\cal D}^{r,t}f=w^{r,t}f$. Therefore, $\rho\circ w^{r,t}f = \rho\circ\s\circ{\cal D}^{r,t}f$. Since $\rho\circ\s$ is 
a $\emptyset$-definable semialgebraic map, this proves the lemma. \ \ \ \ $\Box$\n

{\bf Remark.} Note that in the proof sketch we wrote the partial derivative of a divided difference as a composition
of a semialgebraic function and of $w^{r,t}f$. We do not know a way to write it as a composition of a semialgebraic function 
and of $j^{r,t}f$ (for $f\in C^\infty(\R^n,\R)$ with $n>1$). 
This is the reason why we introduce additional complexity by taking as basic $f$-sets
the preimages of semialgebraic sets under divided differences extensions, 
rather than taking the preimages of semialgebraic sets under multijet extensions, which are simpler.\n

Below, $\pi_{t,m}$,  $t\geq m$, denotes the natural
projection from $D(r,t,n)$ to $D(r,m,n)$ (for which $\pi_{t,m}\circ{\cal D}^{r,t}f={\cal D}^{r,m}f$).\n

{\bf Lemma 2.13.} {\it Let $f\in R_n\subset C^\infty(\R^n,\R)$.
Let $S\subset \R^p \times D(r,m,n)$ be a $\emptyset$-definable semialgebraic set, with a finite $\emptyset$-definable semialgebraic partition $\cal S$.
There exist $t>m$ and a finite $\emptyset$-definable semialgebraic Whitney stratification
${\cal P}$ of $(id\times\pi_{t,m})^{-1}(S)\subset \R^p\times D(r,t,n)$, which refines
$(id \times \pi_{t,m})^{-1}({\cal S})$, such that for each $f\in R_n$:\n
i) the partition $(id\times{\cal D}^{t,m}f)^{-1}({\cal P})$ is a finite Whitney stratification 
which refines the partition $(id\times{\cal D}^{r,m}f)^{-1}({\cal S})$,\n
ii) the upper dimensional strata of $(id\times{\cal D}^{t,m}f)^{-1}({\cal P})$  are monotonic.
}

{\bf Proof.} Let $f\in R_n$. By Proposition 2.11, we may refine $\cal S$ to a finite Whitney stratification ${\cal S}'$ whose
strata are $\emptyset$-definable semialgebraic sets. Since $f\in R_n$, $(id\times{\cal D}^{r,m}f)^{-1}({\cal S}')$ is a finite Whitney stratification.
By Lemma 2.12, the Jacobian of $id\times{\cal D}^{r,m}f$ is a composition of a $\emptyset$-definable semialgebraic map on $\R^p\times D(r,t,n)$ and
of $id\times{\cal D}^{r,t}f$ for some $t>m$. This can be seen to imply that for each coordinate plane $V$ of $\R^p\times(\R^n)^r$ there exists
a $\emptyset$-definable semialgebraic map $\gamma_V:\R^p\times D(r,t,n) \rw \Z$, 
such that the rank of the projection to $V$ of the tangent plane to the preimage of $\cal S$ at point $(v,u)\in\R^p\times(\R^n)^r$, is given
by $\gamma_V\circ(id\times{\cal D}^{r,t}f)$ (if $(v,u)$ is not in the preimage of $S$, we take the rank to be equal to $-1$). 
Since there are only finitely many coordinate planes in $\R^p\times(\R^n)^r$, we may refine $\pi_{t,m}^{-1}({\cal S}')$ into a $\emptyset$-definable 
semialgebraic partition ${\cal Q}$, such that on each element of ${\cal Q}$ the value of $\gamma_V$, for each coordinate plane $V$, is constant. 
By Proposition 2.11, there exists a finite Whitney stratification ${\cal P}$ which refines ${\cal Q}$, whose strata are $\emptyset$-definable 
semialgebraic sets. Since $f\in R_n$, ${\cal A}=(id\times{\cal D}^{r,t}f)^{-1}({\cal P})$ is a Whitney stratification. 

In general, we cannot expect the strata of $\cal A$ to be monotonic. Indeed, let $T\in {\cal P}$ and let $Y\in\pi_{t,m}^{-1}({\cal S}')$ be 
the stratum of ${\cal S}'$ which contains $T$. Take any coordinate plane $V$ in $\R^p\times(\R^n)^r$. Let $A=(id\times{\cal D}^{r,t}f)^{-1}(T)$ and let
$X=(id\times{\cal D}^{r,t}f)^{-1}(Y)$. By construction, the rank of the
projection of $X$ to $V$ is constant on $A$, which however does 
not imply that that the projection of $A$ itself to $V$ is a constant rank map. Nevertheless, this is true 
if $dim(X)=dim(A)$. In particular, this means that the upper dimensional strata of $\cal A$ are monotonic.\ \ \ \  $\Box$\n 

We get to the main point of this section.\n

{\bf Theorem 2.14.} {\it Let $f\in R_n$.
Let $S\subset \R^p \times D(r,m,n)$ be a $\emptyset$-definable semialgebraic set, with a finite $\emptyset$-definable semialgebraic partition $\cal S$.
Then the set $(id\times{\cal D}^{r,m}f)^{-1}(S)$ has a finite Whitney stratification 
which refines $(id\times{\cal D}^{r,m}f)^{-1}({\cal S})$, and whose
strata are monotonic and are $\emptyset$-definable basic $f$-sets in $\R^p\times(\R^n)^r$.
}

{\bf Proof.} Let  $M=(id\times{\cal D}^{r,m}f)^{-1}(S)$.
We make the following induction assumption. 
For each positive integer $i$ there exist integers $t_j$, 
$\emptyset$-definable semialgebraic sets $P_j\subset \R^p \times D(r,t_j,n)$ with finite Whitney stratifications ${\cal P}_j$,
$j=1,..,i$, and a semialgebraic set $S_i\subset\R^p\times D(r,t_i,n)$ with a finite Whitney stratification ${\cal S}_i$, for which the following holds:\n
i) the strata of ${\cal P}_j$ are of equal dimension; denoting their codimension by $c_j$, we have $c_1<c_2<..<c_i<codim({\cal S}_i)$,\n  
ii) the strata of the stratifications ${\cal P}_1,..,{\cal P}_i,{\cal S}_i$ are $\emptyset$-definable semialgebraic sets,\n
iii) the strata (recall that $f\in R_n$) of $(id\times{\cal D}^{r,t_j}f)^{-1}({\cal P}_j)$ are monotonic,\n
iv) the union of $(id\times{\cal D}^{r,t_i}f)^{-1}({\cal S}_i)$ and of $(id\times{\cal D}^{r,t_j}f)^{-1}({\cal P}_j)$, 
$j=1,..,i$ is a Whitney stratification of $M$.\n
Note that the induction claim is true for $i=1$ by Lemma 2.13. Assume that the induction claim is true for $i=N-1$. 
By Lemma 2.13, there exists $t_N>t_{N-1}$ and a finite Whitney stratification ${\cal P}$ of $(id\times \pi^{t_N,t_{N-1}})^{-1}(S_{N-1})$, 
which refines $(id\times \pi^{t_N,t_{N-1}})^{-1}({\cal S}_{N-1})$ and whose strata are $\emptyset$-definable semialgebraic sets,
with the additional property that the upper dimensional strata of ${\cal A}=(id\times{\cal D}^{r,t_N}f)^{-1}({\cal P})$ are
monotonic. Take ${\cal P}_N$ to be the collection of the upper dimensional strata of ${\cal P}$ and take ${\cal S}_N$ to be the 
collection of the rest of strata of $\cal P$. Note that this choice satisfies i),ii),iii) for $i=N$. 
The union $\cal A$ of $(id\times{\cal D}^{r,t_N}f)^{-1}({\cal P}_N)$ and $(id\times{\cal D}^{r,t_N}f)^{-1}({\cal S}_N)$
is a Whitney stratification which refines $ (id\times{\cal D}^{r,t_{N-1}}f)^{-1}({\cal S}_{N-1})$. Since
iv) holds for $i=N-1$ and $dim(S_{N-1})<dim(P_{j})$, $j=1,..,N-1$, Proposition 2.10 implies that iv) holds 
for $i=N$ as well.

Thus the induction claim is true. Since the codimension of $S_i$ grows with $i$, for some $i=I$ its codimension will become larger
than $dim(\R^p\times(\R^n)^r)$. Since $f\in R_n$, this means that $(id\times{\cal D}^{r,t_N}f)^{-1}({\cal S}_I)$ is empty.
Thus the collection  of strata of $(id\times{\cal D}^{r,t_N}f)^{-1}({\cal P}_i)$, $i=1,..,I$, forms, according to iv) and iii),
a Whitney stratification $\cal M$ of $M$ whose strata 
are monotonic. According to ii), each stratum of $\cal M$ is a $\emptyset$-definable basic $f$-set. \ \ \ $\Box$\n

Theorem 2.14 has the following immediate corollary.\n

{\bf Corollary 2.15.} {\it Let $f\in R_n$. Each $\emptyset$-definable basic $f$-set has a finite Whitney stratification whose
strata are monotonic and are $\emptyset$-definable basic $f$-sets.} \ \ \  $\Box$ \n\n\n

{\bf 3. Cylindrical decomposition of projections of basic $f$-sets.}\n

Let $f\in R_n\subset C^\infty(\R^n,\R)$ (for definition of $R_n$, see section 2).
In this section we intend to show that projections of basic $f$-sets on coordinate planes admit cylindrical decomposition whose cells are again 
projections of basic $f$-sets. This then allows to prove Theorem A stated in Introduction. We first establish some auxiliary facts.
Denote by $\pi_{\R^k,P}$ the projection from $\R^k$ to a coordinate plane $P$ in $\R^k$.\n

{\bf Definition 3.1.} {\it 
We call $X\subset\R^k$ a {\it $\emptyset$-definable $f$-set}, if there exists a $\emptyset$-definable basic $f$-set $S\subset\R^p\times(\R^n)^r$,
for some $p,r,n\geq 0$, such that $X$ is the projection of $S$ on a coordinate plane of $\R^p\times(\R^n)^r$.\n
}

{\bf Lemma 3.2.} {\it Let $X,Y\subset\R^k$ be $\emptyset$-definable $f$-sets. Then their intersection and union are $\emptyset$-definable
$f$-sets as well.}

{\bf Proof.} The sets $X,Y$, being $\emptyset$-definable $f$-sets, can be defined by existential formulas of the language ${\cal L}_{\wh{Df}}$.
Therefore, their union and intersection can be also defined by existential formulas of ${\cal L}_{\wh{Df}}$. By Lemma 2.8, both sets
are then the projections on coordinate planes of $\emptyset$-definable basic $f$-sets. \ \ \ \ $\Box$\n 

{\bf Lemma 3.3.} {\it 
Let $X\subset [-1,1]^k$ be a $\emptyset$-definable $f$-set. Then there exists a $\emptyset$-definable
basic $f$-set $S\subset\R^p\times(\R^n)^r$ contained in $[-1,1]^(p+nr)$, for some $p,r,n\geq 0$, 
such that $X$ is the projection of $S$ on a coordinate plane of $\R^p\times(\R^n)^r$.
}

{\bf Proof.} Since $X$ is a $\emptyset$-definable $f$-set, 
there exist a $\emptyset$-definable basic $f$-set $S_0\subset\R^{p'}\times (\R^{n})^{r}$, for some nonnegative integers
$p',r$, and a $k$-dimensional coordinate plane $P$ in $\R^{p'}\times (\R^{n})^{r}$, such that $X$ can be obtained as the projection of $S_0$ to $P$. 
By definition, $S_0$ is the preimage of a semialgebraic set $A\subset\R^{p'}\times D(r,m,n)$ 
under $id\times{\cal D}^{r,m}f:\R^{p'}\times(\R^n)^r\rw \R^{p'}\times D(r,m,n)$, 
for some $m\geq 0$, such that $A\subset\R^{p'}\times [-1,1]^{nr}\times \widetilde{D}(r,m,n)$.
If $S_0\subset[-1,1]^{p'+nr}$, we may take $S=S_0$. If not, choose a coordinate $x_i$ from $\R^{p'}\times(\R^n)^{r}\cong \R^{p'+nr}$, such that    
the projection of $S_0$ on the $x_i$ axis is not contained in $[-1,1]^{p'+nr}$. 
Note that necessarily $1\leq i\leq p'$, since $S_0$ is a subset of $\R^{p'}\times[-1,1]^{nr}$ (because 
$A\subset\R^{p'}\times [-1,1]^{nr}\times \widetilde{D}(r,m,n,)$).
The projection of $\R^{p'+nr}$ to $P$ must be along the coordinate $x_i$, since $X\subset[-1,1]^k$. 
Denote by $V_1$ the coordinate plane in $\R^{p'}\times(\R^n)^r$, obtained by setting $x_i=0$. It is not difficult to check (since $1\leq i\leq p'$) that 
the projection of $S_0$ on $V_1\times(\R^n)^r$ is given by the preimage of the projection of $A$ on $P\times D(r,m,n)$, under 
$id\times{\cal D}^{r,m}:V_1\times(\R^n)^r\rw V_1\times D(r,m,n)$. Since the projection of a $\emptyset$-definable 
semialgebraic set is a $\emptyset$-definable semialgebraic set, we conclude that 
there exists a $\emptyset$-definable basic $f$-set $S_1\subset \R^{p'-1}\times(\R^n)^{r}$, such that $X$ is the projection of $S_1$ on 
the coordinate plane $P$. 

We may repeat this argument, getting a sequence of $\emptyset$-definable basic $f$-sets $S_j\subset V_j\times(\R^s)^r$, $dim(V_j)=p'-j$,
such that for each $j$, $X$ is the projection of $S_j$ on $P$. We may do so until $S_j$ becomes a subset of $[-1,1]^{p'-j+nr}$, which occurs
after $l\leq p'$ steps. Now take $S=S_l$.   \ \ \ \ $\Box$\n   

{\bf Lemma 3.4.} {\it
Let $S\subset\R^p\times(\R^n)^r$ be a $\emptyset$-definable basic $f$-set, contained in $[-1,1]^{p+nr}$. Then there exists a finite monotonic
Whitney stratification $\cal I$ of $[-1,1]^{p+nr}$, whose strata are $\emptyset$-definable basic $f$-sets, 
such that $S$ is stratified by a subset of $\cal I$. 
}

{\bf Proof.} 
By definition, $S$ is the preimage of a $\emptyset$-definable semialgebraic set $A\subset\R^p\times D(r,m,n)$ under $id\times{\cal D}^{r,m}f$, 
for some $m\geq 0$. 
Since $S\subset[-1,1]^{p+nr}$, we may assume that $A$ is contained in $[-1,1]^{p+nr}\times \widetilde{D}(r,m,n)$.
Consider the partition $\cal P$ of $[-1,1]^{p+nr}\times \widetilde{D}(r,m,n)$ whose elements are $A$ and its complement, 
and note that the 
preimage of $[-1,1]^{p+nr}\times \widetilde{D}(r,m,n)$, under $id\times{\cal D}^{r,m}f$, is $[-1,1]^{p+nr}$. 
By Theorem 2.14, there exists 
a finite Whitney stratification $\cal I$ which refines $(id\times{\cal D}^{r,m}f)^{-1}({\cal P})$, and whose strata are 
monotonic and are $\emptyset$-definable basic $f$-sets. Since $\cal I$ refines $(id\times{\cal D}^{r,m}f)^{-1}({\cal P})$,
there is a subset of $\cal I$ which stratifies $S$. \ \ \ \ $\Box$\n

Let $X$ be a $\emptyset$-definable $f$-set. Then there exists 
a $\emptyset$-definable basic $f$-set $S$, such that $X$ is the projection of $S$ on some coordinate plane $P$.  
By Corollary 2.15 $S$ admits a finite stratification with monotonic strata. We define the dimension of $X$ to be equal to the 
maximal rank with which the strata of this stratification project on $P$ (one may check that this definition does not
depend on the choice of $S$ and its stratification).
    
We now state the cylindrical decomposition result.\n 

{\bf Theorem 3.5.} {\it
Let $A\subset \R^s=\R^{s-1}\times\R$ be a $\emptyset$-definable
$f$-set, with the coordinates on $\R^{s-1}\times\R$ being denoted by $(x,t)$. 
Then there exists a partition ${\cal B}$ of $\R^{s-1}$ into finitely many
connected $\emptyset$-definable $f$-sets, such that $\forall B\in{\cal B}$
there is a finite family of continuous functions $g_i:B\rw\R$, $i= 0,1,..,l_B+1$
$$
g_0(x)\equiv-\infty<g_1(x)<...<g_{s_B}(x)<g_{l_B+1}(x)\equiv +\infty,
$$   
with the property that the family of sets of the form
$$
\{(x,t):x\in B, \ g_i(x)<t<g_{i+1}(x)\}
$$
('stripe sets'), and sets of the form
$$
\{(x,t):x\in B, \ t=g_i(x)\}
$$
('graph sets'), constitutes a partition $\cal J$ of $\R^s$ into finitely many connected $\emptyset$-definable 
$f$- sets, and there exists a subset ${\cal J}'$ of $\cal J$ which constitutes a partition of $A$.
}

{\bf Proof.} The proof of the theorem will be by induction on the dimension $s$. In the case $s=1$, the theorem
just says that a $\emptyset$-definable $f$-set $A\subset\R$ consists of finitely many
components. If $A\subset[-1,1]$, then by Lemma 3.3 there exists a $\emptyset$-definable basic $f$-set 
$Q\subset [-1,1]^q\subset \R^p\times(\R^n)^r$, $q=p+nr$, such that $A$ is the projection of $Q$ on a coordinate
plane of $\R^p\times(\R^n)^r$. By Lemma 3.4, there exists a finite Whitney stratification of $[-1,1]^q$,
a subset of which stratifies $Q$. The connected components of the strata form again a finite Whitney stratification.
Therefore $Q$, and thus also its projection, consist of finitely many components. In the case that
$A\not\subset[-1,1]$, $A=(A\cap[-1,1])\cup(A-[-1,1])$. Since by Lemma 3.2 $A\cap[-1,1]$ is a $\emptyset$-definable $f$-set,
we only have to show that $A-[-1,1]$ has finitely many components. Since $A-[-1,1]=A\cap(\R-[-1,1])$, by Lemma 3.2
it is a $\emptyset$-definable $f$-set. It is not difficult to check that the map 
$x\mapsto 1/x$ maps $A-[-1,1]$ into another $\emptyset$-definable $f$-set, which is now contained
in $[-1,1]$. Therefore $A-[-1,1]$, and thus $A$, have finitely many components.\n

We make the induction assumption that the theorem is true in all dimensions
smaller than $s$. Suppose first that $A\subset[-1,1]^{s}$. 
Let us make the following {\it ad hoc} definition.\n

{\bf Definition 3.5.1.} We say that a $\emptyset$-definable $f$-set $G\subset[-1,1]^{s-1}\subset\R^{s-1}$ projects {\it well}   
on a coordinate plane $P$ of $\R^{s-1}$ if the following holds:\n
i) $G$ projects injectively on $P$, and $\pi_{\R^{s-1},P}(G)\subset P$ is open,\n
ii) there exist $p,r\geq 0$, a finite monotonic Whitney stratification ${\cal I}$ of the unit cube $[-1,1]^{q}\subset 
\R^p\times(\R^s)^r$, $q=p+sr$, whose strata are $\emptyset$-definable basic
$f$-sets, and a set $J\subset[-1,1]^q$, stratified by a subset of $\cal I$, such that $\forall x\in G$, 
$$
\pi_{\R^s,\R^{s-1}}^{-1}(x)\cap A \ = \ \pi_{\R^{q},\R^s}\left( \pi_{\R^{q},P}^{-1}\left(  \pi_{\R^{s-1},P}(x) \right) \cap J   \right),
$$
iii) each $y\in\pi_{\R^{s-1},P}(G)$ is a regular value of $\pi_{\R^{q},P}|_{{\cal I}}$.\n

We now make two claims, whose proof we defer until later.\n 

{\bf Claim 3.5.2.} {\it Suppose $A\subset[-1,1]^s$. If the theorem holds in all dimensions smaller than $s$, then
there exists a partition of $\pi_{\R^s,\R^{s-1}}(A)$ into finitely many $\emptyset$-definable $f$-sets
$G_1,..,G_N$, such that for each $G_i$, $i=1,..,N$, there exists a coordinate plane $P_i$ of $\R^{s-1}$
on which $G_i$ projects well.\n
}

{\bf Claim 3.5.3.} {\it  
Suppose that the theorem holds in all dimensions smaller than $s$, and let $G\subset\R^k$, $k<s$, be a $\emptyset$-definable $f$-set. 
Then $G$ has a finite number of components,
and each of them is a \ef. Moreover, if $F_1,..,F_N\subset \R^k$ are \ef s, then $G-F_1\cup ..\cup F_N$ is a \ef.}\n 

We intend to show, assuming Claims 3.5.2 and 3.5.3 
that the induction assumption holds also in dimension $s$. Fix $i$, $1\leq i\leq N$. 
Since $G_i$ projects well on $P_i$, $U=\pi_{\R^{s-1},P_i}(G_i)\subset P_i$ is open,
and there exist $p,r\geq 0$, a finite Whitney stratification ${\cal I}$ of the unit cube $[-1,1]^{q}\subset 
\R^p\times(\R^n)^r$, $q=p+nr$,
whose strata are monotonic $\emptyset$-definable basic
$f$-sets, and a subset ${\cal J}$ of this stratification, such that, writing $J=\cup_{S\in{\cal J}}S$, 
$$
\pi_{\R^s,\R^{s-1}}^{-1}(x)\cap A \ = \ \pi_{\R^{q},\R^s}\left( \pi_{\R^{q},P}^{-1}\left(  \pi_{\R^{s-1},P}(x) \right) \cap J   \right)
$$
holds for each $x\in G_i$.
Moreover, each $y\in U$ is a regular value of $\pi_{\R^{q},P}|_{{\cal I}}$.
Note that $ {\cal I}^{dim \ P_i}_y=\pi_{\R^{q},P_i}^{-1}(y) \ \cap \ {\cal I}^{dim P_i} $ consists of isolated points and is compact.
Thus ${\cal I}^{dim \ P_i}_y$ consists of a finite number of points $\forall y\in U$.
Fix nonnegative integers $j,k_1,..,k_j$. 
Let $U_{j}^{k_1,..,k_j}$ be the subset of points $y\in U$, for which 
the points of $\pi_{\R^q,\R^{s}}\left({\cal I}^{dim \ P_i}_y \right)$ project to precisely $j$ distinct points on $\R^{s-1}$, 
$x_1,..,x_j$, ordered, say, lexicographically,
and the cardinality of $ \pi_{\R^q,\R^{s}}\left({\cal I}^{dim \ P_i}_y \right)\cap\pi_{\R^s,\R^{s-1}}^{-1}(x_l)$ is equal to $k_l$
for each $l=1,..,j$. These sets form a partition of $U$, which we show now to be finite.\n

Indeed, observe that $\pi_{\R^{q},P_i}^{-1}(U)\cap {\cal I}^{dim P_i}$
is a closed Whitney stratified subset of the manifold $\pi_{\R^{q},P_i}^{-1}(U)$, and that   
$\pi_{\R^{q},P_i}$ is a submersion on each stratum $T$ of $\pi_{\R^{q},P_i}^{-1}(U)\cap {\cal I}^{dim P_i}$.
Moreover, the map $\pi_{\R^{q},P_i}|_{\ov{T}}:\ov{T}\rw U$,
where $\ov{T}$ denotes the closure of $T$ in $\pi_{\R^{q},P_i}^{-1}(U)$,
is proper. By the Isotopy Lemma, 
the fibers ${\cal I}^{dim P_i}_y$ are homeomorphic over connected components of $U$. 
Since the fibers are compact and consist of isolated points, there exists $K>0$ such that each
fiber consists of not more than $K$ points. Thus the number of nonempty $U_{j}^{k_1,..,k_j}$ sets is finite. 
The set $U$ is a \ef, so it can be defined by an existential ${\cal L}_{\wh{Df}}$ formula. 
Since the stratification $\cal I$ is finite and consists of $\emptyset$-definable basic $f$-sets,
one may write a suitable ${\cal L}_{\wh{Df}}$ formula for each $U_{j}^{k_1,..,k_j}$ and conclude, by Claim 3.5.3, 
that each $U_{j}^{k_1,..,k_j}\subset P_i$ is a $\emptyset$-definable $f$-set.

By Claim 3.5.3 again, \ $U_j^{k_1,..,k_j}$ has finitely many components
each of which is a $\emptyset$-definable $f$-set. Take any such component and denote it by $H$. Denote by $G_{i,H}\subset\R^{s-1}$ the set
$\pi_{\R^{s-1},P_i}^{-1}(H)\cap G_i$. 
We take $g_l:G_{i,H}\rw\R$, $l=1,..,h$ ($h$ is equal to one of $k_1,..,k_j$) 
to be the functions which send $x\in G_{i,H}$ 
to the projections on the $t$-axis of the points of 
$$ 
\pi_{\R^{q},\R^s}\left({\cal I}^{dim \ P_i}_{y}   \right) \ \cap \ \pi^{-1}_{\R^s,\R^{s-1}}(x),
$$ 
$y=\pi_{\R^{s-1},P_i}(x)$, ordered by magnitude. Observe that the corresponding 'stripe' and 'graph' sets are again $\emptyset$-definable
$f$-sets. 

Since each $y\in H$ is a regular value of $\pi_{\R^q,P_i}|_{\cal I}$,
the set $G_{i,H}$ has the following property, implied by
the Isotopy Lemma. Namely, the points of $\pi_{\R^q,\R^{s}}\left({\cal I}^{dim \ P_i}_y \right)$, $y=\pi_{\R^{s-1},P_i}(x)$,
vary continuously as $x$ varies in $G_{i,H}$. Let $a(x), b(x)$ be two such points, such that 
$\pi_{\R^q,\R^{s-1}}(a(x_0))=\pi_{\R^q,\R^{s-1}}(b(x_0))=x_0$
for some $x_0\in G_{i,H}$. Then it follows from the definition of 
the set $U_j^{k_1,..,k_j}$ and the fact that $H$ is connected, that $\pi_{\R^q,\R^s}(a(x))=\pi_{\R^q,\R^s}(b(x))=x$
for all $x\in G_{i,H}$. In fact, the projections of $a(x)$ and $b(x)$ on the line $\pi_{\R^s,\R^{s-1}}^{-1}(x)$ will be either
equal for all $x\in G_{i,H}$, or distinct for all $x\in G_{i,H}$. This implies in particular that
the functions $g_l$, $l=1,..,h$ are continuous.\n

Denote by ${\cal I}'_y$, $y=\pi_{\R^{s-1},P_i}(x)$, the stratification obtained by taking the components of strata of ${\cal I}_y$.
Note that the set $J_y=\pi_{\R^q,\R^s}^{-1}(y)\cap J$ projects under $\pi_{\R^q,\R^s}$ on $\pi_{\R^s,\R^{s-1}}^{-1}(x)$, and
is stratified by $\pi_{\R^q,\R^s}^{-1}(y)\cap {\cal J}\subset {\cal I}_y$.
Since the frontier condition is satisfied for ${\cal I}'_y$, $\ov{J_y}$ is stratified by a subset of
${\cal I}'_y$, and has therefore a monotonic Whitney stratification with connected strata 
which we denote by $\ov{{\cal J}_y}$. Since $J_y$ is compact and $\pi_{\R^q,\R^s}(\ov{J_y})\subset\pi_{\R^s,\R^{s-1}}^{-1}(x)$,
there exist, by Lemma 1.5, a subset ${\cal T}_x\subset\ov{{\cal J}_y}$ of $1$-dimensional strata, 
which project with rank $1$ to $\pi_{\R^{s},\R^{s-1}}^{-1}(x)$, and a subset ${\cal P}_x\subset\ov{{\cal J}_y}$ of 
$0$-dimensional strata, with the following property. Namely, the sets $T_x=\cup_{S\in {\cal T}_x}\pi_{\R^q,\R^s}(S)$ and
$P_x=\cup_{S\in {\cal P}_x}\pi_{\R^q,\R^s}(S)$ are disjoint and their union is equal to $\pi_{\R^{s},\R^{s-1}}^{-1}(x)\cap A$. 
Moreover, the boundary points of the projections of strata of ${\cal T}_x$ are 
projections of points from ${\cal I}^{dim \ P_i}_y$. 
The Isotopy Lemma implies that the projections of the sets $T_x$ and $P_x$ vary continuously (in the Hausdorff metric)
as $x$ varies in $G_{i,H}$. Together with the fact that
each two continuously varying points from ${\cal I}^{dim \ P_i}_y$, $y=\pi_{\R^{s-1},P_i}(x)$, have projections which are either 
always equal or always distinct, as $x$ varies over $G_{i,H}$, this implies the following: 
if for some $x_0\in G_{i,H}$,  
$(x_0,g_l(x_0))\in \pi_{\R^s,\R^{s-1}}^{-1}(x_0)\cap A$, then 
$(x,g_l(x))\in \pi_{\R^s,\R^{s-1}}^{-1}(x)\cap A$ is true $\forall x\in G_{i,H}$. Similarly, if for some  $x_0\in G_{i,H}$,
$\{x_0\}\times(g_l(x_0),g_{l+1}(x_0))\subset \pi_{\R^s,\R^{s-1}}^{-1}(x_0)\cap A$, then for each  $x\in G_{i,H}$ one has  
$\{x\}\times(g_l(x),g_{l+1}(x))\subset   \pi_{\R^s,\R^{s-1}}^{-1}(x)\cap A$. 
This shows that the partition of $\pi_{\R^s,\R^{s-1}}^{-1}(G_{i,H})$ into 'stripe' and 'graph' sets, generated by the functions $g_1,..,g_l$,
is such that the set $\pi_{\R^s,\R^{s-1}}^{-1}(G_{i,H})\cap A$ is a union of elements from this partition.\n

Thus in the case $A\subset[-1,1]^s$, we may take $\cal B$ to consist of components of
$G_{i,H}$, where $H$ ranges over the components of $U_j^{k_1,..,k_j}$, $U=\pi_{\R^{s-1},P_i}(G_i)$,
$i,j,k_1,..,k_j\in \Z^+$, and of components of the complement $\R^{s-1}-G$. By Claim 3.5.3, 
these are $\emptyset$-definable $f$-sets.
This verifies the induction step in the case $A\subset[-1,1]^s$.\n

Suppose now that $A$ is not a subset of $[-1,1]^s$. Let $A_{i_1,..,i_l}$, $l\leq n$,
$1\leq i_1\leq ..\leq i_l\leq n$, denote the subset of points of $A$ whose $i_1,..,i_l$ coordinates have modulus greater than $1$,
and the rest of their coordinates have modulus equal or less than $1$. These sets form a partition of $A$.
Fix $l\leq n$ and $i_1,..,i_l$, $1\leq i_1\leq ..\leq i_l\leq n$. Let $\a_{i_1,..,i_l}$ denote the mapping which sends
$x_{i_j}$ to $1/x_{i_j}$ for each $j=1,..,l$, and keeps the rest of coordinates unchanged. Note that it maps $A_{i_1,..,i_l}$ to
its homeomorphic image inside $[-1,1]^s$. It is not difficult to see that this image is again a \ef.
Applying to these homeomorphic image the result which we obtained for $A\subset[-1,1]^s$, it is possible, via $\a_{i_1,..,i_l}^{-1}$,
to verify that the theorem is true also for the set $A_{i_1,..,i_l}$ itself. This verifies the induction 
step in the case that $A\not\in[-1,1]^s$, and proves the theorem.\n

It remains to prove Claims 3.5.2 and 3.5.3 on which we relied in the course of the proof of the theorem.\n

{\bf Proof of Claim 3.5.3.} If the theorem holds in dimensions smaller than $s$,
there exists a partition of $\R^{k}$ into finitely many connected \ef s, such that $G$ is a
union of elements from a subset of this partition. Each component of $G$ must be a union of sets from this partition. Hence, by Lemma 3.2,
each component is a \ef. Further, by Lemma 3.2 $F=F_1\cup ..\cup F_N$ is a \ef, and $G-F=(\R^k-F)\cap G$. By the 
induction assumption and Lemma 3.2, $\R^k-F$, being a union of \ef s, is itself a \ef. \ \ \ \ $\Box$\n  

To prove Claim 3.5.2, we first prove an auxiliary statement.\n 

{\bf Claim 3.5.4.} {\it Suppose that $A\subset[-1,1]^s$. 
Let $G\subset\pi_{\R^s,\R^{s-1}}(A)$ be a $\emptyset$-definable $f$-set, and 
let $P$ be a coordinate plane in $\R^{s-1}$, $dim(P)=dim(G)$. If the theorem holds in dimensions smaller than $s$,
there exists a partition of $G$ into finitely
many $\emptyset$-definable $f$-sets $E,G_1,..,G_N$, such that $dim(\pi_{\R^{s-1},P}(E))<dim(P)$,
and each $G_i$, $i=1,..,N$, projects well on $P$.
}

{\bf Proof of Claim 3.5.4.} If $dim(\pi_{\R^{s-1},P}(G))<dim(P)$, then we just take $E=G$. If $dim(\pi_{\R^{s-1},P}(G))=dim(P)$,
then the induction assumption and Lemma 3.2 imply that there exists $K>0$, such that
$\pi_{\R^{s-1},P}^{-1}(y) \cap G$ consists of at most $K$ points $\forall y\in P-E'$, where $E'$ is 
a $\emptyset$-definable set of dimension smaller than $dim(P)$. 
Order the points in $\R^{s-1}$ by, say,
the lexicographical order relation. For each $y\in P-E'$ denote by $x_i(y)$ the $i-th$ largest point of
$\pi_{\R^{s-1},P}^{-1}(y)\cap G$.
Denote by $F_i$, $i=1,..,K$, the set 
$$
\{ x_i(y): \  y\in P-E' \ \ for \ which \ \ |\pi_{\R^{s-1},P}^{-1}(y)\cap G|\geq i\}.
$$  
The sets $F_i$ partition $G\cap\pi_{\R^{s-1},P}^{-1}(P-E')$. It can be seen, applying Claim 3.5.3, 
that the sets $F_i$ are \ef s. Fix $i$, $1\leq i\leq K$. Since $F_i$ is a \ef, there exists a \bef \ $S\subset\R^p \times(\R^n)^r$
and a finite monotonic Whitney stratification ${\cal I}$ of $[-1,1]^q$, $q=p+nr$, such that $S$ is stratified by a subset of $\cal I$, and
$\pi^{-1}_{\R^s,\R^{s-1}}(F_i)\cap A=\pi_{\R^q,\R^s}(S)$. Let $H\subset P$ be the union of projections to $P$ of strata of $S$ which project to $P$
with rank $dim(P)$. By Proposition 1.2 the union of strata of ${\cal I}(dim(P)-1,P)$ is a compact set, and therefore its projection to $P$, which we denote by 
$Z$, is compact. By Lemma 3.2 $Z$ is a \ef. By the induction assumption and Lemma 3.2, $H-Z$ is a \ef. Put $G_i=G\cap\pi_{\R^{s-1},P}^{-1}(H-Z)$.
The set $G_i$ projects well to $P$, and by Lemma 3.2 is a \ef. Note that $\pi_{\R^{s-1},P}(G)-\cup_i{G_i}$ has dimension smaller than $dim(P)$.
Denote the preimage of this set in $G$ by $E$. The sets $G_i$, $i=1,..,K$, and $E$ form a partition of $G$. 
By Claim 3.5.3 \  $E$ is a \ef \ as well. \ \ \ \ $\Box$\n

{\bf Proof of Claim 3.5.2.} We suppose that the dimension of $G$ is $d_1\leq n-1$. Put $E_{11}=G$. Let us enumerate all coordinate planes 
of $\R^{s-1}$ of dimension $d_1$: $P^{d_1}_1,..,P^{d_1}_{k_1}$. There exists a partition of $E_{11}$ into the \ef s $G^{d_1}_{11},..,G^{d_1}_{1N_1}$,
which project well on $P^{d_1}_1$, and a \ef \ $E_{12}$ which projects on $P^{d_1}_1$ with dimension smaller than $d_1$. Repeating the same step with
$E_{12}$ and $P^{d_1}_2$, we obtain $G^{d_1}_{21},..,G^{d_1}_{1N_2}$ and $E_{13}$. Repeating this step $k_1$ times, we obtain a partition of $G$
into the \ef s $G^{d_1}_{ij}$ which project well to $P^{d_1}_i$, $1\leq i\leq k_1$, $1\leq j\leq N_i$, and the \ef \ $E_{1k_1}$, which we 
also denote by $E_{21}$. Note that the dimension of $E_{21}$, which we denote by $d_2$ must be smaller than $d_1$. Enumerate now all
coordinate planes of $\R^{s-1}$ of dimension $d_2$, and repeat for dimension $d_2$ what we have done earlier for dimension $d_1$. 
We continue in this fashion until for some $l\geq 1$ we obtain $E_{l1}=\emptyset$. 
The \ef s $G^{d_i}_{ij}$,
which project well on $P^{d_i}_j$, $1\leq j\leq N_{i-1}$, $i=1,..,l-1$, form a partition of $G$. \ \ \ \ $\Box$\n

This finishes the proof of Theorem 3.5. \ \ \ \ $\Box$\n

We now state our main theorem.\n

{\bf Theorem 3.6.} {\it The theory $T_{\wh{Df}}$ is model complete and o-minimal.}

{\bf Proof.} To show that $T_{\wh{Df}}$ is model complete it is sufficient to show that in the structure $\R_{\wh{Df}}$, 
the complement of any set which is 
a projection of a $\emptyset$-definable quantifier free set, is itself a projection of a $\emptyset$-definable 
quantifier free set. By Lemmas 2.7 and 2.8 this is equivalent to being the complement of any $\emptyset$-definable 
$f$-set again a \ef. The latter is a corollary of Theorem 3.5 and Lemma 3.2. 

Further, to show o-minimality, we have to show that any definable (with parameters) set has finitely many components.
Let such set be denoted by $S$. There is a formula $\phi(y_1,..,y_l)$ 
and $y_{i_1 0},..,y_{i_m 0}\in\R$, such that $S$ is defined by $\phi$ with $y_{i_j}=y_{i_j 0}$, $j=1,..,m$.
Note that $S$ can be identified with the intersection of the set defined by $\phi(y_1,..,y_l)$, 
and the plane given by $y_{i_j}=y_{i_j 0}$, $j=1,..,m$.
By model completeness $\phi(y)$ is equivalent to an existential formula $\exists x \psi(x,y)$, where $\psi(x,y)$ 
is a quantifier free formula. Thus $S$ is the intersection 
of the set $M$ defined by $\exists x \psi(x,y)$, with  the plane given by $y_{i_j}=y_{i_j 0}$, $j=1,..,m$.
Since by Lemma 2.8 $M$ is a \ef, Theorem 3.5 applies, and can be seen to imply the finiteness of the number 
of components of $S$. \ \ \ \ $\Box$\n\n

{\bf 4. A generalization and possible applications.}\n

We comment on how the results extend to the case of at most countably many 
generic smooth functions. This corresponds to taking instead of the map $id\times j^m_rf$ the 
map $id\times j^m_rf_1\times ..\times j^m_rf_k$ and instead of $id \times {\cal D}^{r,m}f$
the map $id\times {\cal D}^{r,m}f_1\times .. \times {\cal D}^{r,m}f_k$. The transversality
arguments go through since the functions (and the divided differences) 
depend on disjoint sets of arguments. Also, the product of $C^{\infty}(\R^{n_i},\R)$, $i=1,2,..$ 
is a Baire space. With these remarks, the proof of Theorem B is almost identical with 
the proof of Theorem A.\n

It seems that Theorem B allows to simplify, at least conceptually, some arguments
which appear in Mather's proof of the topological stability
of proper generic smooth maps ([Ma1], [Ma2], [Ma3]; we refer to the version 
given in [GWPL]). 
It seems that one of the main difficulties in this proof is to show that 
a generic smooth map admits a so called Thom stratification. A related simpler problem 
is to stratify the range of a generic smooth map. Our result gives the following
short proof of the existence of such stratification.\n

{\bf Proposition 4.1.} {\it Let $A\subset[-1,1]^s\subset\R^s$ be a semialgebraic set,
and let $f:\R^s\rw\R^m$ be a generic smooth map. Then $f(A)$ admits a $C^k$
Whitney stratification for any $k\geq 1$.}

{\bf Proof.} We identify $C^\infty(\R^n,\R^m)$ with $(C^{\infty}(\R^n,\R))^m$,
and let $f=(f_1,..,f_m)\in (C^{\infty}(\R^n,\R))^m$. 
The set $f(A)$ is defined by 
$$
\exists x (x\in A\wedge  y=f(x)),
$$
which, since $A\subset[-1,1]^s$, can be easily rewritten as a formula of ${\cal L}_{\wh{Df_1},..,\wh{Df_m}}$.
By the results of the theory of o-minimal structures ([vdDM], [L]), definable
sets in o-minimal structures are $C^k$ Whitney stratifiable, and thus, since $f$ is generic and the conclusion
of Theorem B holds, $f(A)$ is a $C^k$ Whitney stratifiable set. \ \ \ $\Box$\n

We remark that although the proof is short, it in fact relies on the Isotopy Lemma
and on the existence of $C^k$ Whitney stratifications of definable sets in o-minimal structures.\n

A key result about semialgebraic sets is that given two semialgebraic submanifolds $X,Y$, the set of points
of $Y$ at which $X$ is not regular over Y, denoted $B(Y,X)$, is a semialgebraic sets of dimension smaller than $dim(Y)$.
One can show that in fact\n 

{\bf Proposition 4.2.} {\it Let $X,Y$ be $C^k$ submanifolds, $k\geq 1$, definable in $\R_{\wh{Df_1},..,\wh{Df_m}}$,
where $f_1,..,f_m$ are generic smooth functions. Then $B(Y,X)$ is a definable set as well, of dimension smaller than $Y$.}\n

(Note that by Theorem B the dimension is well defined). In [GWPL], Chapter I, section 3, there is a proof that a
generic polynomial map admits a Thom stratification. Proposition 4.2 allows us to repeat this proof 
(with relatively minor modifications), and to conclude, assuming the fact that a generic smooth function
restricted to its critical set is finite to one, that:\n

{\bf Proposition 4.3.} {\it Let $A\subset[-1,1]^s\subset\R^s$ be an open semialgebraic set. Then for a generic smooth $f:\R^s\rw\R^p$,
the map $f|_A:A\rw\R^p$ admits a Thom stratification.}\n

The author got initially interested in the problems discussed in this article during an
attempt to generalize the results of [AgGa] to the generic smooth setting. 
Related questions were raised before in [Suss]. The results presented here are not sufficient 
to answer most of such questions, since one also needs to consider 
(suitably restricted) flows of generic vector fields.\n\n

{\bf Acknowledgements.} This text is a revised and corrected version of the draft which 
the author wrote while staying at SISSA, Trieste. The author thanks  
A. Agrachev for stimulating discussions, and SISSA and Max Planck Institute for Mathematics 
for their financial support. 

\bigskip
\bigskip
\bigskip
\bigskip

\begin{center}{\bf REFERENCES}\end{center}
{\bf [AgGa].} Agrachev, A.; Gauthier, J.-P. On the subanalyticity of Carnot-Caratheodory distances, 
{\it Ann. Inst. H. Poincare Anal. Non Lineaire } 18 (2001), no. 3, 359--382.\n
{\bf [BZ].} Berezin I.S., Zhidkov N.P., {\it Computing Methods}, vol. 1, Pergamon Press, Oxfprd e.a., 1965.\n
{\bf [Ch].} Chatzidakis, Z., Introduction to Model Theory, lecture
notes, Luminy 2001. Available from http://www.logique.jussieu.fr/www.zoe/.\n
{\bf [Gab].} Gabrielov, A. M. Projections of semianalytic sets. (Russian) {\it Funkcional. Anal. i Prilozen.} 2 1968 no. 4, 18--30.\n
{\bf [GG].} Golubitsky, M.; Guillemin, V., {\it Stable mappings and their singularities}, 
Graduate Texts in Mathematics, Vol. 14. Springer-Verlag, New York-Heidelberg, 1973.\n
{\bf [GWPL].} Gibson, C. G.; Wirthmuller, K.; du Plessis, A. A.; Looijenga, E. J. N., {\it Topological stability of smooth mappings}, 
Lecture Notes in Mathematics, Vol. 552. Springer-Verlag, Berlin-New York, 1976.\n 
{\bf [GrY].} Grigoriev, A.; Yakovenko, S., Topology of generic multijet preimages and blow-up via Newton interpolation, {\it 
J. Differential Equations} 150 (1998), no. 2, 349--362. \n
{\bf [vdD1].} van den Dries, L.,  A generalization of the Tarski-Seidenberg theorem, and some nondefinability results, {\it Bull. AMS } 15 (1986), 189-193.\n 
{\bf [vdD2].}  van den Dries, L., O-minimal structures, {\it Logic: from foundations to applications} (Staffordshire, 1993), 137--185, Oxford Sci. Publ., 
Oxford Univ. Press, New York, 1996.\n
{\bf [vdDM].} van den Dries, L.; Miller, C., Geometric categories and o-minimal structures, {\it Duke Math. J.} 84 (1996), no. 2, 497--540.\n
{\bf [IY].} Ilyashenko, Yu.; Yakovenko, S. Finite cyclicity of elementary polycycles in generic families, {\it Concerning the Hilbert 16th problem}, 21--95, 
Amer. Math. Soc. Transl. Ser. 2, 165, Amer. Math. Soc., Providence, RI, 1995.\n  
{\bf [L].} Ta Le Loi, Verdier and strict Thom stratifications in o-minimal structures, {\it Illinois J. Math.} 42 (1998), no. 2, 347--356.\n 
{\bf [KM].} Karpinski, M.; Macintyre, A., A generalization of Wilkie's theorem of the complement, and an application to Pfaffian closure,
{\it Selecta Math. (N.S.)} 5 (1999), no. 4, 507--516.\n
{\bf [Ma1].} Mather, J. N., Stratifications and mappings, {\it Dynamical systems} 
(Proc. Sympos., Univ. Bahia, Salvador, 1971), pp. 195--232. Academic Press, New York, 1973.\n
{\bf [Ma2].} Mather, J. N., How to stratify mappings and jet spaces, 
{\it Singularites d'applications differentiables} 
(Sem., Plans-sur-Bex, 1975), pp. 128--176. Lecture Notes in Math., Vol. 535, Springer, Berlin, 1976.\n 
{\bf [Ma3].} Mather, J. N., {\it mimeographed notes}.\n 
{\bf [RSW].} Rolin, J.-P.; Speissegger, P.; Wilkie, A. J., Quasianalytic Denjoy-Carleman classes and o-minimality, {\it J. Amer. Math. Soc.} 16 (2003), no. 4, 751--777.\n 
{\bf [Suss].}  Sussmann, H.,  Some optimal control applications of real analytic stratification and desingularization,  
{\it Singularities Symposium -- Lojasiewicz 70}, B. Jakubczyk, W.
Pawlucki, and J. Stasica Eds., Banach Center Publications Vol. 44, Polish Academy of Sciences, Warsaw, Poland, 1998, pp. 211-232.\n 
{\bf [W1].} Wilkie, A. J., A theorem of the complement and some new o-minimal structures, {\it Selecta Math. (N.S.)} 5 (1999), no. 4, 397--421.\n
{\bf [W2].}  Wilkie, A. J., Model completeness results for expansions of the ordered field of real numbers by restricted Pfaffian 
functions and the exponential function. {\it J. Amer. Math. Soc.} 9 (1996), no. 4, 1051--1094. \n\n

Present address for correspondence:\n
{\tt alexg@mpim-bonn.mpg.de}

\end{document}